\documentclass[11pt]{article}

\usepackage{epsfig, color, rotating}
\usepackage{amssymb, amsmath, latexsym}

\title{Isometric embeddings of families of special Lagrangian submanifolds.}
\author{Diego Matessi \thanks
{Dipartimento di Scienze e Tecnologie Avanzate, Universit\`{a} del Piemonte
Orientale, Alessandria, Italy
(e-mail: matessi@unipmn.it)}}

\begin{document}

\maketitle

%%%%%%%%%%%%%%%%%%%%%%%%%%%%%%%%%%%%%%%%%%%%%%%%%%%
\newcommand{\ubar}{\ensuremath{\overline{u}}}
\newcommand{\Omr}{\ensuremath{\Omega_{\text{r}}}}
\newcommand{\omr}{\ensuremath{\omega_{\text{r}}}}
\newcommand{\Omst}{\ensuremath{\Omega_{\text{st}}}}
\newcommand{\omst}{\ensuremath{\omega_{\text{st}}}}
\newcommand{\thest}{\ensuremath{\theta_{\text{st}}}}
%%%%%%%%%%%%%%%%%%%%%%%%%%%%%%%%%%%%%%%%%%%%%%%%%%%%%

\newcommand{\pder}[2]
    {\ensuremath{\frac{\partial {#1}}{\partial {#2}}}}
\newcommand{\delvect}[1]
    {\ensuremath{\frac{\partial}{\partial {#1}}}}
\newcommand{\dbar}{\ensuremath{\overline{\partial}}}
\newcommand{\zbar}{\ensuremath{\overline{z}}}
\newcommand{\Ombar}{\ensuremath{\overline{\Omega}}}
\newcommand{\reals}[1]{\ensuremath{\mathbb{R}^{#1}}}
\newcommand{\rationals}[1]{\ensuremath{\mathbb{Q}^{#1}}}
\newcommand{\integ}[1]{\ensuremath{\mathbb{Z}^{#1}}}
\newcommand{\nat}{\ensuremath{\mathbb{N}}}
\newcommand{\complex}[1]{\ensuremath{\mathbb{C}^{#1}}}
\newcommand{\projc}[1]{\ensuremath{\mathbb{C} \mathbb{P}^{#1}}}
\newcommand{\projr}[1]{\ensuremath{\mathbb{R} \mathbb{P}^{#1}}}
\newcommand{\projany}[1]{\ensuremath{\mathbb{P}^{#1}}}
\newcommand{\Img}{\operatorname{Im}}
\newcommand{\Rl}{\operatorname{Re}}
\newcommand{\Vl}{\operatorname{Vol}}
\newcommand{\Rc}{\operatorname{Ric}}
\newcommand{\Hes}{\operatorname{Hess}}
\newcommand{\Id}{\operatorname{Id}}
\newcommand{\tr}{\operatorname{Tr}}
\newcommand{\gr}[1]{\operatorname{Gr}_{#1}}
\newcommand{\inner}[2]{ \langle {#1}, {#2} \rangle}
\newcommand{\gener}{\operatorname{span}}
\newcommand{\Endo}{\operatorname{End}}
\newcommand{\So}{\operatorname{SO}}
\newcommand{\Su}{\operatorname{SU}}
\newcommand{\su}{\mathfrak{su}}
\newcommand{\so}{\mathfrak{so}}
\newcommand{\Gl}{\operatorname{GL}}
\newtheorem{defin}{Definition}
\newtheorem{lem}{Lemma}
\newtheorem{thm}{Theorem}
\newtheorem{cor}{Corollary}
\newtheorem{prop}{Proposition}
\newtheorem{conjecture}{Conjecture}
\newenvironment{exmpls}{\vspace{2ex} \hspace{-\parindent}
                            \textbf{Examples:}
                            \begin{enumerate}}{\end{enumerate}}
\newenvironment{proof}[1]{ \textbf{#1}}{\hfill $\Box$\vspace{2ex}

                                                     }
\newenvironment{rmk}[1]
       {\vspace{2ex}\hspace{-\parindent}
          \textbf{Remark {#1}.}}{\hfill $\Box$\vspace{2ex}

                                                      }
\newenvironment{observation}[1]
 {\vspace{2ex}\hspace{-\parindent}
       \textbf{{#1}}}{\hfill $\Box$ 
                              \vspace{2ex}
                             
                           }

%comands specific to this paper
\newcommand{\Chom}{\hat{\Omega}}
\newcommand{\chom}{\hat{\omega}}
\newcommand{\cod}{\operatorname{codim}}
\begin{abstract}
We prove that certain Riemannian manifolds can be isometrically
embedded inside Calabi-Yau manifolds.
For example we prove that given any real-analytic one parameter 
family of Riemannian metrics $g_t$ on a $3$-dimensional manifold $Y$ 
with volume form independent of $t$ and with a real-analytic family
of nowhere vanishing harmonic one forms $\theta_t$, then $(Y, g_t)$
can be realized as a family of special Lagrangian
submanifolds of a Calabi-Yau manifold $X$. 
We also prove that certain principal torus bundles can be 
equivariantly and isometrically embedded inside Calabi-Yau manifolds
with torus action. We use this to construct examples
of $n$-parameter families of special Lagrangian tori inside 
$n+k$-dimensional Calabi-Yau manifolds with torus symmetry.
We also compute McLean's metric of $3$-dimensional 
special Lagrangian fibrations with $T^2$-symmetry.
\end{abstract}

\section{Introduction}
Since the work of McLean it is known that compact
special Lagrangian submanifolds with positive first Betti number $b_1$
always come in families, i.e. can always be deformed to form
a $b_1$-dimensional moduli space. So for example it may happen that a
a family of special Lagrangian tori locally fibres the Calabi-Yau manifold.
After the SYZ conjecture in mirror symmetry (cfr.  
\cite{syz:slfib, hitch:msslag, mgross:slfibtop, mgross:slfibgeom}) it is
believed that the geometry of the ambient Calabi-Yau could be understood 
in terms of the families of special Lagrangian tori that it contains. 
In this paper we find ways to construct some examples of families of 
sL submanifolds. For instance we show that a large class of $1$-parameter or 
$2$-parameter families of Riemannian $3$-tori can be realised as families 
of special Lagrangian tori inside a $3$-dimensional Calabi-Yau manifold. 
These results generalise those of Bryant and of the author.
The former proved in \cite{bry:calembed} that any real-analytic Riemannian 
$3$-dimensional manifold can be isometrically embedded in a Calabi-Yau 
manifold as a special Lagrangian manifold. In \cite{diego:slonepar} we 
refined this result to show that this could be done also for some examples of
real-analytic $1$-parameter families of Riemannian tori. 

A Calabi-Yau manifold (CY for short) is a K\"{a}hler $n$-dimensional manifold 
$X$, with a nowhere vanishing holomorphic $n$-form $\Omega$ and a Ricci-flat
K\"{a}hler form $\omega$. We will refer to the pair $(\Omega, \omega)$ as
a CY structure. When $\omega$ is not Ricci-flat then we will call it an 
almost Calabi-Yau structure (ACY). We will find often convenient to use
the following characterisation of Calabi-Yau manifolds due to Hitchin:

\begin{thm} [Hitchin] \label{cyforms}
If $X$ is a $2n$-dimensional oriented manifold and $\Omega$ a
complex valued $C^{\infty}$ $n$-form satisfying:
\newcounter{uno}
\begin{list}{(\roman{uno})}{\usecounter{uno} \setlength{\parsep}{0cm} 
\setlength{\topsep}{\itemsep} \setlength{\leftmargin}{.5cm} }
 \item  $d\Omega = 0$; 
 \item locally there exist $1$-forms $\theta_{1}, \ldots, \theta_{n}$
such that $\Omega =\theta_{1} \wedge \ldots \wedge \theta_{n}$
($\Omega$ is said to be locally decomposable); 
 \item $|\Omega \wedge \Ombar|> 0 $;
\end{list}
then $\Omega$ determines an integrable almost complex structure
on $X$ with respect to which $\Omega$ is a holomorphic
$(n,0)$-form.

Moreover, suppose $\omega$ is a symplectic form on $X$ such that
\begin{list}{(\roman{uno})}{
\usecounter{uno} 
\setlength{\parsep}{0cm} 
\setlength{\topsep}{\itemsep} 
\setlength{\leftmargin}{.5cm} 
                                }
\setcounter{uno}{3}  
\item $\Omega \wedge \omega = 0$;
 \item  $\omega^{n} = c \, \Omega \wedge \Ombar$, for some constant $c$;
\end{list}
then $\omega$ is a Ricci-flat K\"{a}hler form with respect to the 
complex structure induced by $\Omega$.
\end{thm}
When $(v)$ is satisfied, we will in general assume that $\Omega$ has been 
rescaled so that the following normalised relation holds
\[ \omega^{n}/n! = 
          (-1)^{n(n-1)/2}(i/2)^{n} \, \Omega \wedge \Ombar. \] 
If only condition (i)-(iv) are satisfied then $(\Omega, \omega)$ is an 
ACY structure.

An $n$-dimensional real submanifold $Y$ of a CY manifold $X$ is special 
Lagrangian (sL for short) if
\[ \Rl \Omega|_{Y} = \Vl_{Y}, \]
where $\Vl_{Y}$ denotes the volume form on $Y$. Equivalently, $Y$ is 
sL if 
\begin{eqnarray}
  \Img \Omega |_{Y} & = & 0 \\
         \omega|_Y  & = & 0.
\end{eqnarray} 
In the case of an ACY manifold we will always use the latter as the definition
of sL submanifold.

We now describe some important results on the geometry of sL submanifolds
to which we will often refer in this paper. In \cite{mclean:deform} McLean
proved that given a compact sL submanifold $Y$ of a CY manifold $X$, the 
moduli space $M_Y$ of nearby sL submanifolds can be identified with a
smooth finite dimensional submanifold of $\Omega^1(Y)$, the space of
$1$-forms of $Y$. The tangent space to $M_Y$ is the space $\mathcal{H}^1(Y)$
of harmonic $1$-forms on $Y$, therefore $M_Y$ has dimension $b_1(Y)$. If 
$Y_t$ is a variation of $Y$ through sL submanifolds such that $Y_0 = Y$ and
$V$ is a normal vector field along $Y$ tangent to this variation then
\[ \theta_V = \iota_V \omega |_Y \]
is a harmonic form and it is tangent to $M_Y$. We call it the variational harmonic
form. McLean proved that every harmonic form gives rise to such a variation
$Y_t$. He also observed that $M_Y$ has a natural metric. If $\theta_1$ and
$\theta_2$ are two harmonic forms then
\[ \inner{ \theta_1}{\theta_2} = \int_Y \theta_1 \wedge \star \theta_2, \]
the standard $L^2$ inner product on $\mathcal{H}(Y)$.

Given $q \in M_Y$ denote $Y_q$ the sL submanifold corresponding to $q$.
If $M_Y$ is connected then we can identify $H_1(Y_q, \integ{})$ with
$H_1(Y, \integ{})$. Now suppose $\dim M_Y = m$ and we are given a smooth 
choice of a basis of harmonic $1$-forms $\theta_1(q), \ldots, \theta_{m}(q)$ on
$Y_q$. Fix a basis $\Sigma_1, \ldots, \Sigma_m$ for the free part of 
$H_1(Y, \integ{})$ and define
\[ P_{ij} (q) = \int_{\Sigma_i} \theta_j. \]
The matrix $P = (P_{ij})$ is called the period matrix. 

Suppose now that $\theta_1(q), \ldots, \theta_{m}(q)$ are chosen so
that $P = \Id$, then we can define the function
\begin{eqnarray}
 \Phi & : & M_y \rightarrow \reals{}  \nonumber \\
  \   & \ & q \mapsto \det( \inner{\theta_i(q)}{\theta_j(q)}). \label{semf:vol}
\end{eqnarray}
This function was introduced by Hitchin in \cite{hitch:msslag}. 
He showed that $T^{\ast}M_Y$ possesses a natural ACY structure. 
He also proved that this ACY structure is CY if and only if $\Phi$ is constant. 
In the case of a family of $2$ dimensional sL tori in a $2$ dimensional CY 
manifold, he also proved in \cite{hitch:msclag} that  $\Phi$ is in fact constant. 
In \cite{diego:slonepar} we gave a $3$-dimensional example where this 
function is not constant. We will call $\Phi$ the semi-flat volumes function.

\vspace{2ex}

\textbf{Acknowledgements.} Early versions of the results in this article 
were contained in my PhD thesis written at Warwick 
University under the supervision of
Mario Micallef and Mark Gross from 1998 to 2001. 
I wish to thank them for their invaluable help and time. 
I have worked on the new versions in this article
sporadically throughout my post-docs, first with an ISM-CRM grant
at Centre de Recherche Math\'{e}matiques, Montr\'{e}al, Canada and then with
an EPSRC grant at Imperial College, London, UK. The paper was completed 
while I was a guest at the Mathematics Departments of the University of
Turin and Politecnico, Italy. I wish to thank these institutions for the 
hospitality. I also thank: Sung Ho Wang for the help in understanding 
the Cartan-K\"{a}hler Theorem and in proving Theorem~\ref{onep:sl} 
and Vestislav Apostolov for some useful discussions. 

\section{Summary of main results}
In this paper we prove various existence results of CY structures
which also imply that some families of Riemannian manifolds can
be embedded isometrically inside CY manifolds. 
 In \cite{bry:calembed} Bryant proved that if $(Y,g)$ is a real-analytic
$3$-dimensional manifold with a real-analytic Riemannian metric $g$ then
there exists a $3$-dimensional CY manifold $X$ and an isometric embedding 
$\iota: Y \rightarrow X$ such that $\iota(Y)$ is a sL submanifold. 
The proof is an application of the Cartan-K\"{a}hler theorem. 
Here we use the same method to extend the result to $1$-parameter families.
In Theorem~\ref{onep:sl} we prove (again in the $3$-dimensional case)
that if $g_t$ is a real-analytic $1$-parameter
family of metrics on $Y$ with constant volume form admitting a real-analytic
family of harmonic forms $\theta_t$ with no zeroes, then there exists a 
CY manifold $X$ and a smooth family of isometric embeddings 
$\iota_t: Y \rightarrow X$ such that $\iota_t(Y)$ is sL for all $t$. 
Moreover $\theta_t$ it the variational harmonic form associated to 
the variation $\iota_t$. The same result was proved by the author
in \cite{diego:slonepar} in the special case of tori and for a more
restricted choice of variational harmonic form.
The theorem is proved in Section~\ref{onp}.
We review the Cartan-K\"{a}hler theorem and Bryant's theorem in 
Sections~\ref{ckt} and \ref{brt}.

In Section~\ref{cy:actins} we study torus actions on an ACY manifold
which preserve the structure $(\Omega, \omega)$. We reformulate the
Pedersen and Poon's ansatz \cite{pp:KEconstr} for K\"{a}hler-Einstein 
metrics with isometric and Hamiltonian torus actions in the case of 
structure preserving actions.  In Section~\ref{isbun} we use it to prove
isometric embedding theorems for torus bundles. We prove that certain
principal torus bundles $E \rightarrow N$ with connection $\Theta$ over
an ACY manifold $N$ can be equivariantly and isometrically embedded in 
a higher dimensional Calabi-Yau manifold $X$ with structure preserving 
torus action. 
In Section~\ref{s1:exmpls}, using these results and a theorem of Goldstein 
\cite{goldstein:cftorus} 
and Gross \cite{mgross:slexmp} on the lifting of sL fibrations, 
we construct some examples of families of Riemannian manifolds which
can be isometrically embedded in CY manifolds. 

In the last Section we assume that we have a solution to the
ansatz in the case of 2-torus symmetry on a $3$-dimensional
CY manifold $X$ and a sL torus fibration $f:X \rightarrow B$ with
fibres containing the orbits of the action. For this fibration
we compute the periods, McLean's metric and the semi-flat volumes
function. We obtain interesting formulae. Using the isometric 
embedding result for $2$-torus bundles we generalise the example
in \cite{diego:slonepar} showing that the semi-flat volume function
may not be constant. 

\section{The Cartan-K\"{a}hler theorem} \label{ckt}
An exterior differential system is a pair $(M, \mathcal{I})$
where $M$ is a smooth $n$-dimensional manifold  and $\mathcal{I}$ 
an ideal of its (graded)-ring of smooth differential forms $\Omega^{\ast}(M)$
which is closed under exterior differentiation. We will 
often denote the exterior differential system just by $\mathcal{I}$
when the underlying manifold $M$ is understood.
We denote $\mathcal{I}_k = \mathcal{I} \cap \Omega^{k}(M)$.
An integral submanifold of $(M, \mathcal{I})$
is a submanifold $\iota: N \rightarrow M$ such that 
$\iota^{\ast} \phi = 0$ for every $\phi \in \mathcal{I}$. 
Many problems in geometry can be stated in terms of finding 
integral submanifolds of an exterior differential system. 

Let $p \in M$, a $k$-dimensional subspace $E$ of $T_pM$ is 
called an integral element of $\mathcal{I}$ if $\phi_{|E} = 0$ 
for all $\phi \in \mathcal{I}$. Let 
$V_{k}(\mathcal{I}) \subseteq \gr{k}(TM)$ be the set of all 
$k$-dimensional integral elements of $\mathcal{I}$. A priori 
$V_{k}(\mathcal{I})$ has no obvious nice structure, for
example it may very well be the whole Grassmannian at one point
and empty at all the others. There are however notions of 
$V_k(\mathcal{I})$ being well behaved near a point in
$\gr{k}(TM)$. 
Let $E \in \gr{k}(T_pM)$ be an integral element of $\mathcal{I}$,
we can choose coordinates $(x_j)$ on a neighbourhood 
$U$ of $p$ in $M$ so that 
$(dx_1 \wedge \ldots \wedge dx_k)_{|E} \neq 0$.
There is an open neighbourhood $W$ of $E$ in $\gr{k}(TM)$
such that for any $E^{\prime} \in W$, 
$(dx_1 \wedge \ldots \wedge dx_k)_{|E^{\prime}} \neq 0$, 
therefore to any $\phi \in \mathcal{I}_k$, we can associate
a map $F_{\phi}: W \rightarrow \reals{}$ such that
\[ \phi_{|E^{\prime}} 
      = F_{\phi}(E^{\prime}) dx_1 \wedge \ldots \wedge dx_k. \]
We have that $ E^{\prime} \in W$ is integral if and
only if $F_{\phi}(E^{\prime}) = 0$ for all $\phi \in \mathcal{I}_k$.
\begin{defin}
We say that $E$ is an \textbf{ordinary element}
if there are forms 
$\phi_1, \ldots, \phi_l$ in $\mathcal{I}_k$ and an open subset
$W^{\prime} \subseteq W$ so that $F_{\phi_1}, \ldots, F_{\phi_l}$
have linearly independent differentials on $W^{\prime}$ and
\[ W^{\prime} \cap V_k(\mathcal{I}) = 
      \{ E^{\prime} \in W^{\prime} \, | \,
      F_{\phi_1}(E^{\prime}) = \ldots = F_{\phi_l}(E^{\prime}) = 0 \}. \]
\end{defin}
In particular if $E$ is an ordinary element then  $V_k(\mathcal{I})$
is a submanifold of $\gr{k}(TM)$ near $E$.

Given an integral element $E \in \gr{k}(T_pM)$ we would like
to measure the extent to which $E$ can be enlarged to a 
higher dimensional integral element. The following definition
is quite natural

\begin{defin} The \textbf{polar space} of $E \in \gr{k}(T_pM)$ is the subspace 
of $T_pM$ defined by
\[ H(E) = \{ v \in T_pM \, | \, (\iota_v \phi)_{|E} = 0, \ 
      \text{for all} \ \phi \in \mathcal{I}_{k+1} \}. \]
It is the union of all integral elements containing $E$.   
\end{defin}
Clearly any $k+1$-dimensional integral element of $\mathcal{I}$
containing $E \in \gr{k}(T_pM)$ can be written as  $E + \reals{}v$
with $v \in H(E)$, but $v \notin E$. Therefore the set of $k+1$-dimensional 
integral elements of $\mathcal{I}$ containing $E$ is in one to one
correspondence with $\mathbb{P}( H(E) / E)$. 
\begin{defin} The \textbf{extension rank} of $E$ is
\[ r(E) = \dim (\mathbb{P}( H(E) / E)) = \dim (H(E)) - k -1 \].
\end{defin}
If $r(E) = -1$, then we cannot extend $E$ to a bigger
integral element. 
\begin{defin} An ordinary integral element $E$ is
said to be \textbf{regular} if $r$ is a constant function
on a small neighbourhood of $E$ in $V_k(\mathcal{I})$.
\end{defin}
We can now state
\begin{thm}[The Cartan-K\"{a}hler Theorem]
Let $(M,\mathcal{I})$ be a 
real-ana\-ly\-tic exterior differential system.
Let $N$ be a connected, real-analytic $k$-di\-men\-sio\-nal regular 
integral submanifold of $M$, with non-negative extension rank 
$r$. Let $Z$ be a real-analytic submanifold of co-dimension $r$
in $M$ containing $N$ and such that $\dim(T_pZ \cap H(T_pN)) = k + 1$
at every $p \in N$.

Then there exists a $k+1$-dimensional, real-ana\-ly\-tic integral
manifold $Y$, such that $N \subset Y \subset Z$. Moreover
any other real-analytic integral submanifold $Y^{\prime}$
such that $N \subset Y^{\prime} \subset Z$ coincides with
$Y$ in a neighbourhood of $N$. 
\end{thm}
Although it may be difficult to prove that an integral
submanifold is regular, Cartan devised a powerful criterion
which reduces the problem to linear algebra. We will now 
describe it. 

A $q+1$-tuple of integral elements $F = (E_0, \ldots, E_q)$ of
$\mathcal{I}$ satisfying 
\[ 0 = E_0 \subset E_1 \subset \ldots \subset E_q \subset T_pM \]
and $\dim E_i = i$ is called an \textbf{integral flag}.
If in addition all $E_i$'s are regular up to $E_{q-1}$
it is called a \textbf{regular flag}. Let  $c_{i}(F)$ be the co-dimension of 
$H(E_i)$ in $T_pM$ and $C(F) = \sum_{i = 0}^{q-1} c_i(F)$.
\begin{thm} [Cartan's Test] Let $F = (E_0, \ldots, E_q)$
be an integral flag. Then $V_{q}(\mathcal{I})$ lies in a 
co-dimension $C(F)$ submanifold of $\gr{q}(TM)$ near $E_q$. 
Moreover $V_{q}(\mathcal{I})$ is a submanifold of codimension 
$C$ near $E_q$ if and only if $F$ is a regular flag. 
\end{thm}

\section{SL isometric embedding of manifolds} \label{brt}
We will now describe how the problem of constructing local
Calabi-Yau structures can be interpreted as one
of finding integral submanifolds of a certain exterior 
differential system $\mathcal{I}$ on a manifold.
We will then show how Bryant used this setting to prove his
theorem on isometric embeddings of special Lagrangian 
submanifolds using the Cartan-K\"{a}hler theorem. 

Let $\Omega_0$ and $\omega_0$ 
be the standard holomorphic and K\"{a}hler form respectively
on $ \complex{n} = \reals{2n}$, i.e. the forms
\[ \Omega_0 = dz_1 \wedge \ldots \wedge dz_n, \]
\[ \omega_0 = \frac{i}{2} \sum_{j} dz_j \wedge d\zbar_j. \]
The subgroup of $\Gl(2n, \reals{})$
which preserves $\Omega_0$ and $\omega_0$ is $\Su(n)$, moreover
the sub-ring $(\Lambda^{\ast} \reals{2n})^{\Su(n)}$ of 
$\Lambda^{\ast}\reals{2n}$  of forms
which are preserved by $\Su(n)$ is precisely the one generated
by $\Omega_0$ and $\omega_0$. 

Suppose now we are given a $2n$-dimensional oriented Riemannian 
manifold $(M,g)$. Let $\pi: F \rightarrow M$ 
denote the $\So(2n)$ co-frame bundle 
of $M$, i.e. the principal $\So(2n)$-bundle whose fibre over a
point $p \in M$ consists of all the oriented linear isometries
$u: T_pM \rightarrow \reals{2n}$. We say that the manifold $M$
admits a $\Su(n)$ structure if there exists a principal 
$\Su(n)$ subbundle $P$ of $F$. This is the same as saying
that there exists a covering $\{ U_{\alpha} \}$ of $M$ with open 
subsets and for every $\alpha$ a local co-frame 
$e^{\alpha} =(e^{\alpha}_{1}, \ldots, e^{\alpha}_{2n})$
on $U_{\alpha}$ such that on overlaps $U_{\alpha} \cap U_{\beta}$, 
$e^{\alpha}$ is related to $e^{\beta}$ by a $\Su(n)$ 
transformation. An alternative way to interpret this is the following.
The group $\Su(n)$ acts on $F$. Let us
denote by $S$ the manifold $F / \Su(n)$. 
A choice of covering with $\Su(n)$-coframes gives a section 
$\sigma: M \rightarrow S$, and vice versa a section $\sigma$ gives 
such a covering and therefore a $\Su(n)$-subbundle $P_{\sigma}$.

On $F$ there exists a natural $\reals{2n}$ valued $1$-form $\theta$
defined by
\[ \theta_u(v) = u(\pi_{\ast} v) \]
for every $u \in F$ and $v \in T_uF$. We can use this form to define
forms $\hat{\Omega}_0$ and $\hat{\omega}_0$ on $F$ by
\[ \Chom_0 (v_1, \ldots , v_n) = \Omega_0 
          ( \theta(v_1), \ldots, \theta (v_n)), \]
\[ \chom_0 (v, w) = \omega_0(\theta(v), \theta(w)). \]
The forms $\Chom_0$ and $\chom_0$ descend to the
manifold $S$ and if $M$ admits a $\Su(n)$ structure, we can use the
section $\sigma: M \rightarrow S$ defining the structure to pull
back these forms to genuine forms on $M$ 
\[ \Omega_{\sigma} = \sigma^{\ast} \Chom_0, \]
\[ \omega_{\sigma} = \sigma^{\ast} \chom_0. \]
Which shows that $M$ admits a $\Su(n)$ structure if and 
only if it admits forms $(\Omega, \omega)$ (not necessarily closed)
satisfying (ii) (iii) and (iv) in Theorem~\ref{cyforms}. If in addition we want 
$(\Omega_{\sigma},\omega_{\sigma})$ to define a Calabi-Yau
structure we must require $\sigma$ to be such that 
\[ \sigma^{\ast} d\Chom_0 = 0, \]
\[ \sigma^{\ast} d\chom_0 = 0, \]
i.e. we require the manifold $\sigma(M)$ to be an integral
$2n$-dimensional submanifold of the exterior differential system
$(S, \mathcal{I})$ where $\mathcal{I}$ is generated by $d \Chom_0$ 
and $d\chom_0$.

With this setting we can use the Cartan-K\"{a}hler Theorem to prove
existence of local Calabi-Yau structures, but in order to do so 
effectively we also need to understand the Cartan Test in this 
contest. We may regard the ideal $\mathcal{I}$
as living on $F$ as well as on $S$. Let $\tau: F \rightarrow S$
and  $\overline{\pi}: S \rightarrow M$ be the projections
and denote the bundle $\ker \tau_{\ast}$ by $\mathbf{h}$. The
fibres of $\mathbf{h}$ are canonically isomorphic to $\su(n)$ .
We are interested in integral elements $E \in \gr{k}(T_uF)$ 
(or $E \in \gr{k}(T_{\tau(u)}S)$) which are mapped isomorphically 
onto $T_{\pi(u)}M$ by $\pi_{\ast}$ (or by $\overline{\pi}_{\ast}$). 
Let us denote the set of such elements by $V_{k}(\mathcal{I}, \pi)$
(or $V_{k}(\mathcal{I}, \overline{\pi})$ respectively). 
Given $E \in \gr{k}(TS)$ we denote 
by $\tilde{E}$ a lift of $E$, i.e. an element in $\gr{k}(TF)$ 
such that $\tau_{\ast}( \tilde{E}) = E$. Clearly $\tilde{E}_{1}$, $\tilde{E}_{2}$
are lifts of the same element $E$ if and only if 
$\tilde{E}_{1} + \mathbf{h} = \tilde{E}_{2} +  \mathbf{h}$. 
One can verify that the following holds
\[ H( \tilde{E}) = \tau_{\ast}^{-1}(H(E)). \]
We deduce that 
\begin{equation} \label{cs}
   c(\tilde{E}) = \cod_{T_uF}(H(\tilde{E})) = 
               \cod_{T_{\tau(u)}S}(H(E)) = c(E), 
\end{equation}
which is useful because it is often simpler to compute these 
numbers on lifts.

 Now suppose $E \in V_{2n}( \mathcal{I}, \overline{\pi})$
with a lift $\tilde{E} \in V_{2n}( \mathcal{I}, \pi)$.
There is a natural integral flag 
$\tilde{B} = (\tilde{E}_0, \ldots, \tilde{E}_{2n})$
with $\tilde{E}_{2n} = \tilde{E}$. 
In fact let $\epsilon_1, \ldots, \epsilon_{2n}$
be the dual standard basis in $\reals{2n}$, then we can define
\[ \tilde{E}_k = \{ v \in \tilde{E} \, | \, 
        \epsilon_{k+1} ( \theta(v)) = \ldots 
                = \epsilon_{2n}(\theta(v)) = 0 \} . \]
Clearly $\tau_{\ast}$ maps $\tilde{B}$ to a well
defined flag $B = ( E_0, \ldots, E_{2n})$ with $E_{2n} =E$.
We want to apply the Cartan Test to prove that $B$ is regular.
First of all one can show that near $E$ and $\tilde{E}$,
$V_{2n}( \mathcal{I}, \pi)$ and
$V_{2n}( \mathcal{I}, \overline{\pi})$  are submanifolds of 
$\gr{2n}(TF)$ and $\gr{2n}(TS)$ respectively of the same co-dimension
$2n q$, where $q$ is the co-dimension of $\Su(n)$ in $\So(2n)$.
For example, when $n=3$, which is the case we are interested in,
they are submanifolds of co-dimension $42$. 
To show that $B$ is regular we have to compute the numbers 
$c_k = \cod (H(E_k))$ and show that $C(B) = \sum_{k=0}^{2n-1}c_k$
is equal to $2nq$. By (\ref{cs}) one can do this computation 
for $\tilde{B}$ , which is easier. 
Clearly $\tilde{E} \subseteq H(\tilde{E}_k)$. 
Therefore $H(\tilde{E}_k) = \tilde{E} + \mathbf{h}_k$, where 
$\mathbf{h}_k$ is a 
subspace of the vertical bundle, i.e. a suspace of $\so(2n)$. 
For example suppose
$x \in \so(2n)$ and $v_1, v_2, \ldots, v_{n} \in \tilde{E}_k$, then
\begin{eqnarray*}
 d \Chom (x,v_1, \ldots, v_{n}) & = &
         \mathcal{L}_{x} \Chom (v_1, \ldots , v_{n}) \\
 & = & \mathcal{L}_{x} \Omega_0 ( \theta(v_1), \ldots, \theta(v_n)),
\end{eqnarray*}
where $\mathcal{L}$ denotes the usual Lie derivative. Similar
equalities clearly hold for $d \chom$. If we now think of 
$\reals{k}$ as standardly embedded in $\reals{2n}$, i.e. 
as the set of vectors whose last $2n-k$-coordinates are zero, then
we see that
\[ \mathbf{h}_k = \{ x \in \so(n) \, | \, 
            (\mathcal{L}_x \Omega_0)_{| \reals{k}} = 
            (\mathcal{L}_x \omega_0)_{| \reals{k}} = 0 \}. \]
Therefore
\[ c_k = \cod_{\so(2n)} ( \mathbf{h}_k). \]
So, computing these numbers is a completely algebraic problem. 
In \cite{bry:calembed}, Bryant points out that whether 
$C(B) = 2nq$ holds or not is a property both of $\Su(n)$ and of
the choice of its representation. 
For example he checks that indeed $C(F) = 42$ in the three
dimensional case when we identify $\complex{n}$
with $\reals{2n}$ via $z \mapsto (x,y)$, where $z = x + i y$, i.e.
when the complex structure $J$ on $\reals{2n}$ is the 
matrix
\[ J = \left( \begin{array}{cc}
                 0 & I_n \\
                 -I_n & 0 
              \end{array} \right). \]
In this case one says that $\Su(n)$ is \textbf{regularly presented}.
If one represents $\Su(n)$ in a
different way then it may not satisfy
$C(B) = 2nq$ (see \cite{bry:calembed} for more details).

We are now ready to sketch Bryant's theorem.
\begin{thm} [Bryant, \cite{bry:calembed} ] \label{isoemb:bry}
Let $(Y,g)$ be a $3$-dimensional oriented, connected
real-analytic manifold
with real-analytic Riemannian metric $g$. Then there exists a three
dimensional Calabi-Yau
manifold $X$ and an isometric embedding $\iota: Y \rightarrow X$
such that $\iota(Y)$ is special Lagrangian. 
\end{thm}
\begin{proof}{Sketch of proof.}
It is known that every $3$-dimensional Riemannian manifold is 
parallelisable, moreover, when $M$ and $g$ are
real-analytic, one can
find a real-analytic orthonormal parallelisation. This means
we can assume that we have a globally defined, real-analytic, oriented
frame $e = (e_1, e_2, e_3)$ on $Y$. When $Y$ is embedded in some
Calabi-Yau manifold as a special Lagrangian submanifold, 
$(Je_1, Je_2, Je_3)$ is an orthonormal frame for
the normal bundle $\nu(Y)$ of $Y$. Via the exponential map we
know that a neighbourhood of $Y$ in $X$ is diffeomorphic
to a neighbourhood of $Y \times \{0 \}$ in $M = Y \times \reals{3}$.
Moreover, if 
$(\epsilon_1, \epsilon_2, \epsilon_3, \phi_1, \phi_2, \phi_3)$ is
the co-frame dual to $(e_1, e_2, e_3, Je_1, Je_2, Je_3)$, then 
on $\iota^{\ast} TX$, where $\iota: Y \rightarrow X$ is the embedding,
we have
\[ \Omega = (\epsilon_1 + i \phi_1) \wedge (\epsilon_2 + i \phi_2)
               \wedge (\epsilon_3 + i \phi_3), \]
\[ \omega = \sum_{k=1}^{3} \epsilon_k \wedge \phi_k .\]

Following these observations, it is reasonable to set the 
problem in the following way.
Let $M= Y \times \reals{3}$ and let $(\phi_1, \phi_2, \phi_3)$
be the canonical co-frame on $\reals{3}$, then with the choice of 
global co-frame $(\epsilon_1, \epsilon_2, \epsilon_3, 
\phi_1, \phi_2, \phi_3)$ on $M$ we let $F = M \times \Gl_6(\reals{})$
and $S = M \times ( \Gl_{6}(\reals{}) / \Su(3))$.
We also define the forms $\Chom_0$ and $\chom_0$ on $F$ (and on $S$) 
to be
\[ \Chom_0 = (\epsilon_1 + i \phi_1) \wedge (\epsilon_2 + i \phi_2)
               \wedge (\epsilon_3 + i \phi_3), \]
\[ \chom_0 = \sum_{k=1}^{3} \epsilon_k \wedge \phi_k .\]
Let $\sigma_0: Y \times \{0 \} \rightarrow S$ be defined by 
$\sigma_0(y, 0) = \tau((y,0), I_6)$, where $I_6$ is the identity matrix.
Denote $\sigma_0(Y \times \{ 0 \} )$ by $Y_0$. 
The goal is to find an open neighbourhood $X$ of 
$Y \times \{0 \}$ in $M$ and a section $\sigma: X \rightarrow S$ 
such that $\sigma(y,0) = \sigma_0(y,0)$ and such that $\sigma(X)$
is an integral submanifold of the exterior differential system
generated by $d \Chom_0$ and $d \chom_0$ on $S$. Clearly
$\sigma_{0}^{\ast} \chom_0 = 0$ and 
$\sigma_{0}^{\ast}(\Chom_0) = \epsilon_1 \wedge \epsilon_2 
                        \wedge \epsilon_3$,
and therefore $Y_0$ is an integral manifold.

With a bit of thought one can see that at each $p \in Y_0$
the tangent space $T_pY_0$ is contained in some 
$E \in V_6( \mathcal{I}, \overline{\pi})$ and it is the $E_3$ in the
canonical flag ending in $E$. In particular $Y_0$ is a 
regular submanifold. 

Bryant computes explicitly the (linear)-equations defining
$\mathbf{h}_3$ and finds that it has dimension $31$ in 
$M_6(\reals{})$, i.e. that $c_3 = 5$ and that $H(T_pY_0)$ has
dimension $29$ (i.e. $\dim \mathbf{h}_3 + \dim M - \dim \Su(3)$).
Therefore the extension rank of $Y_0$ is $25$. 
We want to apply the Cartan-K\"{a}hler Theorem to extend $Y_0$
to a $4$-dimensional integral manifold whose tangent spaces
are in $V_4(\mathcal{I}, \overline{\pi})$. We need to find 
a $9$ dimensional submanifold $Z_0$ of $S$, containing $Y_0$,
such that $T_p Z_0 \cap H(T_p Y_0)$ is of dimension $4$ and 
is in $V_4(\mathcal{I}, \overline{\pi})$.

The idea to construct $Z_0$ is the following. Being a linear
subspace of $M_6(\reals{})$, $\mathbf{h}_3$ has a complementary
subspace $W_0$, i.e. a subspace of  $M_6(\reals{})$ of dimension
$5$ such that $W_0 \cap \mathbf{h}_3 = \{ 0 \} $. There exists
a neighbourhood $U_0$ of $0$ in $W_0$ such that for all $x \in U_0$, 
$I_6 + x$ is in $\Gl_n(\reals{})$, therefore we can define
the following submanifold of $F$:
\[ \tilde{Z}_0 = 
   \{ (p,(y_1, 0,0), I_6 + x) \in 
             Y \times \reals{3} \times \Gl_6(\reals{}) \ | \
                   p \in Y, x \in U_0 \}. \]
It is a submanifold of dimension $9$, and since 
$\su(3) \subset \mathbf{h}_{3}$, $W_0$ is transversal
to the $\Su(3)$ orbit, therefore $\tilde{Z}_0$ maps down
to a well defined $9$ dimensional submanifold $Z_0$ of $S$, 
containing $Y_0$. One can check that indeed 
$T_p Z_0 \cap H(T_p Y_0)$ is of dimension $4$ and it is in 
$V_4(\mathcal{I}, \overline{\pi})$.
In fact $T_p Z_0 \cap H(T_p Y_0)$ is sent via $\overline{\pi}_\ast$
isomorphically onto the subspace of $T_{(p,0)}M$ spanned by 
$T_pY$ and $\partial / \partial y_1$. We can now apply 
the Cartan-K\"{a}hler Theorem to find an integral, $4$ dimensional
submanifold $Y_1$, such that $Y_0 \subset Y_1 \subset Z_0$
and whose tangent spaces are all contained in 
$V_4(\mathcal{I}, \overline{\pi})$. In particular we may assume
that $Y_1$ is the image of a section $\sigma_1$ of $S$ defined
over $\{ (p,(y_1, 0,0)) \in Y \times \reals{3} \ | \ y_1 \ \text{is 
sufficiently small} \}$.

We now repeat similar arguments to extend $Y_1$ to a $5$ dimensional
integral submanifold containing $Y_1$ which is the image
of a section of $S$ defined over a $5$ dimensional submanifold in 
$M$. One observes that each tangent space of $Y_1$ is contained
in some $E \in V_6( \mathcal{I}, \overline{\pi})$ and that it
is the $E_4$ of the canonical flag of $E$, therefore $Y_1$ is
a regular submanifold. One computes the equations defining
$\mathbf{h}_4$ and finds that it is a $24$ dimensional subspace
of $M_6(\reals{})$. Therefore $H(T_pY_1)$ is $22$ dimensional and
$T_pY_1$ has extension rank $17$. We need to find $Z_1$, a submanifold
of $S$, of dimension $17$ containing $Y_1$ and with the
property that $T_p Z_1 \cap H(T_p Y_1)$ is of dimension $5$ at every
$p$ and it is contained in $V_5(\mathcal{I}, \overline{\pi})$. 
One finds a subspace $W_1$ of $M_6(\reals{})$ of dimension $12$ 
which is complementary to $\mathbf{h}_4$ and a neighbourhood $U_1$ of
$0$ in $W_1$ such that $I_6 + x \in \Gl_6(\reals{})$ for all $x \in U_1$. 
Then we can define
\[ \tilde{Z}_1 = 
   \{ ((p,y_1, y_2,0), I_6 + x) \in 
             Y \times \reals{3} \times \Gl_6(\reals{}) \ | \
                   p \in Y, x \in U_1 \}. \]
We have that
$\tilde{Z}_1 $ maps down to a well defined $17$ dimensional submanifold
$Z_1$ in $S$, containing $Z_0$ and therefore $Y_1$ which has the
desired properties. We apply Cartan-K\"{a}hler and we find a
five dimensional integral submanifold $Y_2$ containing $Y_1$ 
which is the image of a section $\sigma_2$ defined over the set 
$\{ (p,y_1, y_2,0) \in Y \times \reals{3} \ | \ y_1, y_2\ \text{are 
sufficiently small} \}$. 

Finally one can extend one more time to find the desired 
six dimensional integral submanfold $Y_3$ which is the image
of a section $\sigma$ defined over an open neighbourhood $X$
of $Y \times \{ 0 \}$ in $M$. This can be done since $h_5$ is 
$14$ dimensional (i.e. $c_5 = 22$), implying $H(T_p Y_2)$ is 
$12$ dimensional and the extension rank of $Y_2$ is $6$. 
The arguments above carry through. 
The forms $\sigma^{\ast} \Chom_0$ and $\sigma^{\ast} \chom_0$
define a Calabi-Yau structure on $X$ with respect to which
$Y \times \{ 0 \}$ is special Lagrangian and induces the
given metric $g$. 
\end{proof}

\section{SL isometric embedding of $1$-parameter families} \label{onp}
We now show how Bryant's result can be refined to prove
an isometric embedding theorem of one parameter families
of Riemannian manifolds as one parameter families of 
special Lagrangian submanifolds in a Calabi-Yau manifold.

First let us recall some important facts. Suppose that
$Y$ is a $3$ dimensional manifold and suppose that 
for some $\kappa > 0$ there exists a map 
$I: Y \times [0, \kappa) \rightarrow X$ into some 
$3$-dimensional Calabi-Yau manifold $(X, \Omega, \omega)$,
such that for each $t \in [0, \kappa)$ the map
$I_t = I( \cdot, t)$ is an immersion and $Y_t = I_t(Y)$
is special Lagrangian. We call $I$ a special Lagrangian
variation of Y. The vector field in $X$ along $Y_t$ defined
by
\[ V_t(p) = \pder{I}{t} (p,t), \]
is called the variational vector field of $I$. If $V_t$ is 
normal to $Y_t$ for all $t$ we call $I$ a normal
(special Lagrangian) variation. Let
\[ \theta_t = I_{t}^{\ast}( \iota_{V_t} \omega ). \]
MacLean has shown that $\theta_t$ is a harmonic 
one-form on $Y$ with respect to the metric $g_t = I_{t}^{\ast}g$
where $g$ is the metric on $X$. We call $\theta_t$ the 
variational harmonic form of the variation $I$. 

We now prove
\begin{thm} \label{onep:sl}
Let $Y$ be a real-analytic $3$-dimensional oriented manifold.
Suppose we are given the data of a pair 
$(g_t, \theta_t)_{ t \in (-\kappa, \kappa)}$
where $g_t$ is a real-analytic one parameter family of metrics 
on $Y$ and $\theta_t$ a real-analytic one parameter family of
$1$-forms, satisfying
\newcounter{wuno}
\begin{list} {(\roman{wuno})}{\usecounter{wuno} 
\setlength{\parsep}{0cm} 
\setlength{\topsep}{\itemsep} 
\setlength{\leftmargin}{.5cm}}
\item  $\Vl_{g_t}$, the volume form with respect
to $g_t$, is independent of $t$, i.e.
\[ \pder{}{t} \Vl_{g_t} = 0; \]
\item $\theta_t$ is harmonic with respect to $g_t$ for all 
              $t \in (-\kappa, \kappa)$;
\item $\theta_t (p) \neq 0$, for all 
              $t \in (-\kappa, \kappa)$ and $p \in Y$.
\end{list}
Then there exists a $3$-dimensional 
Calabi-Yau manifold $(X,\Omega, \omega)$
and a normal special Lagrangian variation
$I: Y \times (-\kappa, \kappa) \rightarrow X$
such that for every $t \in (-\kappa, \kappa)$, 
the map $I_t = I( \cdot, t)$ is an isometry
with respect to $g_t$ and $\theta_t$ is the variational harmonic 
one-form with respect to $I$. 
\end{thm}
\begin{proof}{Proof.} Let $X_1 = Y \times (-\kappa, \kappa)$. 
Define the following one forms on $X_1$:
\[ \epsilon_1 (p, t) = \frac{\theta_{t}}{|\theta_{t}|}, \]
\[ \phi_1 (p, t) = | \theta_{t} | dt, \]
where the length of $\theta_{t}$ is computed w.r.t. $g_t$ and
it is non-zero because of (iii).
We can complete these two forms to a global, real-analytic
co-frame $(\epsilon_1, \epsilon_2, \epsilon_3, \phi_1)$ such
that for every $t \in (-\kappa, \kappa)$,
 $(\epsilon_1(\cdot, t) , \epsilon_2 ( \cdot, t), 
                                  \epsilon_3( \cdot, t))$
is an oriented, orthonormal co-frame on $Y$ with respect to 
$g_t$.

Now define $M = X_1 \times \reals{2}$ and let $(y_1, y_2)$ be
coordinates on $\reals{2}$. Let $(\phi_2, \phi_3) = (dy_1, dy_2)$
be the standard co-frame on $\reals{2}$.
Let the co-frame bundle $F$ of $M$ be trivialised by
$(\epsilon_1, \epsilon_2, \epsilon_3, \phi_1, \phi_2, \phi_3)$,
so that we can identify $F$ with $M \times \Gl_6(\reals{})$ and 
$S$ with $M \times (\Gl_6(\reals{}) / \Su(n))$.
As usual, define on $F$ (and on $S$) the forms
\[ \Chom_0 = (\epsilon_1 + i \phi_1) \wedge (\epsilon_2 + i \phi_2)
               \wedge (\epsilon_3 + i \phi_3), \]
\[ \chom_0 = \sum_{k=1}^{3} \epsilon_k \wedge \phi_k .\]
We let $\mathcal{I}$ be the exterior differential system 
on $F$ (or $S$) generated by $d \Chom_0$ and $d \chom_0$.

Now consider the section of $S$ defined on $X_1 \times \{ 0 \}$
by $\sigma_1 ((p,t), 0) = \tau((p,t), 0, I_6)$. We show
that $Y_1 = \sigma_1( X_1 \times \{ 0 \})$ is an integral
submanifold of $(S, \mathcal{I})$. In fact we have
\[ \sigma_{1}^{\ast}( \chom_0) = 
              \epsilon_1 \wedge \phi_1 = \theta_t \wedge dt \]
and
\begin{eqnarray*}
 \sigma_{1}^{\ast}( \Chom_0) 
        & = & \epsilon_1 \wedge \epsilon_2 \wedge \epsilon_3 +
                         i \phi_1 \wedge \epsilon_2 \wedge \epsilon _3 \\
        & = & \epsilon_1 \wedge \epsilon_2 \wedge \epsilon_3 +
                 i | \theta_t| dt \wedge \epsilon_2 \wedge \epsilon _3
\end{eqnarray*}
By assumption (ii), $\theta_t$ is closed on $Y$ for every $t$, 
therefore $\sigma_{1}^{\ast}( \chom_0)$ is also closed on $X_1$. 
Now observe that $\epsilon_1 \wedge \epsilon_2 \wedge \epsilon_3$
restricts to the volume form w.r.t. $g_t$ on 
$Y \times \{ t \} \subset X_1$, therefore, by assumption (i), it is 
independent of $t$, i.e. it is closed on $X_1$. 
Moreover,
if $\star$ is the Hodge-star operator w.r.t the metric $g_t$,
then
\[  \star \theta_t = |\theta_t| \epsilon_2 \wedge \epsilon_3, \]
therefore, by assumption (ii), $ |\theta_t| \epsilon_2 \wedge \epsilon_3$ 
is closed on $Y$. It follows
that $| \theta_t| dt \wedge \epsilon_2 \wedge \epsilon _3$, and 
therefore $\sigma_{1}^{\ast}( \Chom_0)$, is closed in $X_1$. 
We conclude that $Y_1$ is an integral submanifold. 

To complete the proof we can now follow exactly the steps in 
Theorem \ref{isoemb:bry}, starting from the second extension.
One replaces the $Y_1$ in Theorem \ref{isoemb:bry}, obtained
by extending $Y_0$, with the $Y_1$ we have defined here and
follows the argument word by word from then on. 
\end{proof}

\section{Structure preserving torus actions} \label{cy:actins}
The action of a Lie group $G$ on an $n$-dimensional
ACY manifold $(X, \Omega, \omega)$ is called structure preserving if it 
is Hamiltonian and preserves $\Omega$. 
Clearly structure preserving actions 
induce holomorphic isometries on $X$, although the converse
is not in general true as we may have holomorphic actions
which do not preserve $\Omega$. 

We restrict to the case when $G$ is the $m$-torus $T = U(1)^{m}$,
with $m < n$, and the action is free.
Denote by $\mathfrak{t} = i \, \reals{m}$ the 
Lie algebra of $T$. Let $\eta_{1}, \ldots, \eta_{m} \in \Gamma(TX)$ 
be the vector fields induced by the action corresponding to 
the standard basis $ie_{1}, \ldots, ie_{m}$ of $\mathfrak{t}$.
We denote by $\mu_{j} : X \rightarrow \reals{}$ a Hamiltonian
of $\eta_{j}$, i.e. a function such that
 $\iota_{\eta_{j}} \, \omega = d\mu_{j}$.
The function $\mu: X \rightarrow \reals{m}$ given by 
$\mu = (\mu_{1}, \ldots, \mu_{m})$ is the Hamiltonian of the action.
For any $t = (t_{1}, \ldots, t_{m}) \in \mu(X)$ we can form the manifold 
$N = \mu^{-1}(t)/ T^{m}$ 
together with its symplectic form $\omr$ given by
symplectic reduction. 
On $N$ we can also naturally define an $(n-m)$-form $\Omr$
as follows. Denote by $\pi: \mu^{-1}(t) \rightarrow N$
the projection. Let $v_{1}, \ldots, v_{n-m} \in T \mu^{-1}(t)$.
Then
\[ \Omr(\pi_{\ast}v_{1}, \ldots, \pi_{\ast} v_{n-m})
 = \iota_{\eta_{1}} \ldots \iota_{\eta_{m}} \Omega(v_{1}, \ldots, v_{n-m}). \]
 
Goldstein \cite{goldstein:cftorus} and Gross~\cite{mgross:slexmp}
proved the following theorem:
\begin{thm} [Goldstein, Gross] \label{s1:gold}
The manifold $(N, \Omr, \omr)$ is almost Calabi-Yau. Moreover, suppose 
there exists a $T^{m}$-invariant continuous map $g: X \rightarrow B$ to 
an $(n-m)$-dimensional real manifold $B$ such that the induced maps 
$g: N \rightarrow B$ are sL fibrations for all (or almost all) 
$t \in \mu(X)$. Then $(g,\mu): X \rightarrow B \times \reals{m}$ is a sL 
fibration on $X$.
\end{thm}
This theorem is our main motivation to study CY manifolds
with $T$ action. 

In \cite{pp:KEconstr} Pedersen and Poon prove an
ansatz for K\"{a}hler-Einstein metrics with a torus acting through
holomorphic isometries. Here we wish to formulate the ansatz for the 
stronger case of structure preserving actions on CY manifolds.
We start by describing the data required in the ansatz.

Let $\pi: E \rightarrow N$ be a principal $T$-bundle over the
complex $(n-m)$-dimensional manifold $N$. We identify $\mathfrak{t}$
with $i \reals{m}$. Let
$U \subseteq \reals{m}$ be a connected open subset with coordinates
$t = (t_1, \ldots, t_m)$.
On $E$ we have the vertical vector fields $\eta_{1}, \ldots, \eta_{m}$
associated to the standard basis of $i \reals{m}$. 
Let $W = (w_{jk})$ be an $m \times m$ positive definite
symmetric matrix of smooth functions on $N \times U$. 
Let $\omr(t)$ be a smooth family of K\"{a}hler forms on $N$
parametrised by $t \in U$ and 
$\theta(t) \in \Omega^1(E, i\reals{n})$ a smooth family of connection
forms, written in components as $\theta = (\theta_1, \ldots, \theta_m)$. 
Assume also that there exists on $N$ a nowhere vanishing holomorphic 
$(n-m)$-form $\Omr$ such that $( \Omr, \omr(t))$ is an ACY structure 
for all $t$. On the $2n$-dimensional real manifold $X = E \times U$ 
we can define the following $T$-invariant forms
\begin{eqnarray}
& \omega =\pi^{\ast} \omr-i \, \theta_{j} \wedge dt_{j}, \label{tm:om} \\
& \Omega = (i)^m \left( \bigwedge_{j=1}^{m} ( w^{jk} dt_k - \theta_{j}) \right) 
                \wedge \pi^{\ast} \Omr. \label{tm:Om}
\end{eqnarray}
\begin{prop} [Pedersen-Poon,\cite{pp:KEconstr}]\label{t2:forms}
Let $\pi: E \rightarrow N$ be a principal $T$-bundle over the
complex $(n-m)$-dimensional manifold $N$ and let
$W = (w_{jk})$, $(\Omr, \omr)$, 
and $\theta$ be the data defined above. 
Then the pair of forms $(\Omega, \omega)$
on $X = E \times U$ written in (\ref{tm:om}) and (\ref{tm:Om})
form an almost Calabi-Yau structure on $X$ if and only if the data
satisfy
\begin{eqnarray}
      & \frac{\partial}{\partial t_{j}}\omr =  i \, d\theta_{j},
                                      \label{t2:eq1} \\
      & \  \nonumber \\
  &\left[ \frac{\partial \, \theta_i}{\partial t_{j}} 
          \right]_{(0,1)} =
                 - \overline{\partial} \, w^{ij} \label{t2:eq6} \\ 
      & \ \nonumber \\ 
   & \frac{\partial w^{ij}}{\partial t_{k}}
         -\frac{\partial w^{ik}}{\partial t_{j}} =0, \label{t2:eq4} 
\end{eqnarray}
where in the second equation, $[ \cdot ]_{(0,1)}$ is the projection onto
the $(0,1)$ part with respect to the complex structure on $N$. 
The operators $d$ and $\dbar$ are computed on $N$.

The structure is CY if $(\Omr, \omr)$ satisfy
\begin{equation} \label{omr:Omr}
 \omr^{n-m} / (n-m)! =  (-1)^{(n-m)(n-m-1)/2}(i/2)^{n-m} \, \det W^{-1} 
   \, \Omr \wedge \overline{\Omega}_{\text{r}}. 
\end{equation}
Moreover any ACY manifold $(X, \Omega, \omega)$ with a free structure preserving 
$T$-action can be described locally by such a construction.
\end{prop}
\begin{proof}{Proof.} We use Theorem~\ref{cyforms}.
Clearly $\Omega$ is locally decomposable
since this is true of $\Omr$ . Notice that since $\theta$ has values
in $i \reals{m}$
\[ \overline{\theta} = - \theta. \]
It is not difficult to prove the identity
\begin{equation} \label{ugh}
 \Omega \wedge \Ombar = (-2)^{m} \ \det W^{-1} \ 
     dt_1 \wedge \ldots \wedge dt_m 
                        \wedge \theta_1 \wedge \ldots \wedge \theta_m 
\wedge \pi^{\ast}(\Omr \wedge \overline{\Omega}_{\text{r}}).
\end{equation}
Since $\theta_1 \wedge \ldots \wedge \theta_m$ restricts to a volume
form on the vertical space of $E$, (\ref{ugh}) implies condition
(iii) of Theorem~\ref{cyforms}. Condition (iv) easily follows from
the identity
\[  \theta_j \wedge dt_j = - \frac{w_{jl}}{2} 
       (w^{lk} dt_k - \theta_{l})
                        \wedge (w^{jm} dt_m + \theta_{j}). \]
A tedious but straight forward computation, which uses (\ref{ugh}),
also shows that if $\omr$ and $\Omr$ satisfy (\ref{omr:Omr}) then
$\omega$ and $\Omega$ satisfy condition $(5)$.
Finally one can easily check that $\Omega$ and $\omega$ are closed
if and only if equations (\ref{t2:eq1}), (\ref{t2:eq6}) and (\ref{t2:eq4})
are satisfied. 

We now prove that every ACY manifold with a structure preserving
free $T$ action can be described by the ansatz.
Let $\mu = (\mu_1, \ldots, \mu_n)$ be the moment map of the
action. Then, for $t \in \mu(X)$ let $E_t = \mu^{-1}(t)$ and
$N_t = \mu^{-1}(t) / T$ with the reduced ACY structure $(\Omr(t), \omr(t))$.
Let $J$ denote the almost complex structure of $X$. Define
\begin{equation} \label{w}
  w_{jk} = \inner{\eta_{j}}{\eta_{k}} = \omega(\eta_{j},J\eta_{k}), 
\end{equation} 
then for fixed $t$, the matrix $W = (w_{jk})$ is a symmetric 
positive definite matrix of functions on $N_t$. 
Now denote 
\begin{equation} \label{xi}
 \xi_{j} = w^{jk}J\eta_{k},
\end{equation} 
The $\xi_j$'s are linearly independent vector fields orthogonal to 
$\mu^{-1}(t)$. Using the fact that the $\eta_j$'s are Killing holomorphic
fields one can show that
\[ [\xi_j, \xi_k] = 0. \]
Therefore the flows $\Phi_{\xi_{j}}^{t_{j}}$ of the $\xi_{j}$'s, 
which exist for sufficiently small $t_j$'s, define coordinates 
$t= (t_1, \ldots, t_m)$ on the leafs of the distribution defined by
the $\xi_j$'s. Let
\[ \Phi^{t} = \Phi_{\xi_{1}}^{t_1} \circ 
         \ldots \circ \Phi_{\xi_{m}}^{t_m}. \] 
The fact that 
\begin{eqnarray*}
  d \mu_{i}( \xi_{j}) & = & \omega(\eta_{i},  w^{jk} \, J\eta_{k}) \nonumber \\
                   & = & w_{ik}w^{jk} = \delta_{ij}.  \label{dtx} 
\end{eqnarray*}
implies that 
\[ \mu(\Phi^{t} (p) )= t \]
when $p \in \mu^{-1}(0)$. Therefore $\Phi^{t}$ identifies $E_0 = \mu^{-1}(0)$
 with $E_t = \mu^{-1}(t)$. Now, since $\eta_k$ commutes with $\xi_j$,
$\Phi^{t}$ is also equivariant, i.e. $N_0$ and $N_t$ are also naturally
identified. Let $E= E_0$ and $N = N_0$. If we take $U \subseteq \reals{m}$
to be an open neighbourhood of $0$ such that $\Phi^{t}$ is defined for 
$t \in U$, then $\Phi^{t}$ identifies $E \times U$ with $\mu^{-1}(U)$
in an equivariant way, moreover $\mu$ is identified with the projection 
onto $U$. 

The reduced ACY structure on $N_t$ can be considered as a family of
ACY structures $(\Omr(t), \omr(t))$ on $N$. We now prove
\begin{equation} \label{Om:const}
 \frac{\partial}{\partial t_j} \Omr = 0.
\end{equation}
Goldstein \cite[Lemma 1]{goldstein:cftorus} proved that 
$\iota_{\eta_1} \ldots \iota_{\eta_m} \Omega$ is closed.
Therefore (\ref{Om:const}) follows from
\[ \mathcal{L}_{\xi_{j}}(\iota_{\eta_1} \ldots \iota_{\eta_m} \, \Omega) =
     \iota_{\xi_{j}} \,  d \, \iota_{\eta_1} \ldots \iota_{\eta_m} \Omega +
  d \, \iota_{\xi_{j}} \iota_{\eta_{1}} \ldots \iota_{\eta_m} \, \Omega = 0. \]
where the second summand vanishes because $\Omega$ is of type $(n,0)$.

The Riemannian metric on $X$ also induces a family of connections 
$\theta(t)= (\theta_1(t), \ldots, \theta_m(t))$ on $E$. 
Proving that $\Phi^{t}$ identifies $\Omega$ and $\omega$ on $X$
with the forms in (\ref{tm:Om}) and (\ref{tm:om}) respectively
is just a matter of a simple verification. The fact that $(\Omr, \omr)$
always satisfy (\ref{omr:Omr}) whenever $(\Omega, \omega)$ 
satisfy property (v) in Theorem~\ref{cyforms} is also
a rather straight forward computation.  
\end{proof}

\begin{rmk}{1} The way we stated the Pedersen-Poon ansatz is different
from how it was stated in the original paper. This is partially due to 
the fact that we are only considering the case of Ricci-flat metrics
where the action preserves also the holomorphic form. Our formulation
is also more convenient for the applications in the following sections.
\end{rmk}

\section{Isometric embeddings of torus bundles} \label{isbun}
We now prove, using essentially Proposition~\ref{t2:forms} and
the Cauchy-Kowalesky theorem, the following isometric embedding
theorem for principal $T^m$ bundles over an ACY base.
\begin{thm} \label{cy:exist}
Let $\pi: E \rightarrow N$ be a principal $T^m$-bundle over an
$(n-m)$-dimensional ACY manifold $(N, \Omega_{0}, \omega_{0})$ 
with connection $\theta_{0}$ whose curvature is of type $(1,1)$. 
Let $W_0$ be a positive definite $m \times m$ symmetric matrix 
of functions on $N$ such that
\begin{equation} \label{omo:Omo}
 \omega_{0}^{n-m} / (n-m)! =  (-1)^{(n-m)(n-m-1)/2}(i/2)^{n-m} \, 
  \det W^{-1}_{0} \, \Omega_0 \wedge \overline{\Omega}_{0}. 
\end{equation}
If all the data is real-analytic,
then there exists a neighbourhood $U \subset \reals{m}$ of $0$
and a $T^m$-invariant CY structure $(\Omega, \omega)$
on $E \times U$, inducing the connection $\theta_{0}$ and the matrix
$W_0$ on $E \times \{ 0 \}$. Moreover the projection 
onto $U$ is the Hamiltonian of the action and the reduced space 
over $t=0$ is $(N, \Omega_{0}, \omega_{0})$.
\end{thm}
\begin{proof}{Proof.}
By Proposition~\ref{t2:forms} we can assume $\Omr = \Omega_0$.
For given $\omr$, denote by $f_{\omr}$ the unique function such
that
\begin{equation*}
 \omr^{n-m} / (n-m)! =  (-1)^{(n-m)(n-m-1)/2}(i/2)^{n-m} \, 
  f_{\omr} \, \Omr \wedge \Ombar_{r}. 
\end{equation*}
For any $m$-tuple of positive numbers $\epsilon = 
(\epsilon_1, \ldots, \epsilon_m)$ denote 
\[ U_l = \{ (x_1, \ldots, x_m) \in \reals{m} \ | \ 
        x_j < \epsilon_j, \ 1 \leq j \leq l;  \ x_j = 0, 
                           \ l+1 \leq j \leq m \}. \]
For fixed $1 \leq l \leq m$ consider the following system
of equations
\begin{equation*}
(S_l) \ \ \ \left\{ \begin{array}{l}
             \frac{\partial}{\partial t_{l}} \omr =  i \, d\theta_{l},
                                      \\
  \left[ \frac{\partial \, \theta_i}{\partial t_{l}} 
          \right]_{(0,1)} =
                 - \overline{\partial} \, w^{il}, \ \ i=1, \ldots, m \\ 
    \frac{\partial w^{ij}}{\partial t_{l}}
         = \frac{\partial w^{il}}{\partial t_{j}}, 
                                  \ \ j \leq i, \ j \leq l \\ 
    \frac{\partial w^{ij}}{\partial t_{l}}
         = \frac{\partial w^{jl}}{\partial t_{i}}, \ \ i \leq j \leq l \\
        \det W^{-1} = f_{\omr}
              \end{array} \right.                
\end{equation*} 
We regard this system as an evolution equation with time variable
$t_l$. It is a well defined system on $N \times \reals{l}$. We
now show that by subsequently applying the Cauchy-Kowalesky Theorem, 
we can find an $m$-tuple $\epsilon = (\epsilon_1, \ldots, \epsilon_m)$ 
such that for every $l=1, \ldots, m$, there is a solution 
$( \omr, W, \theta)_l$ to $(S_l)$ on $U_l$. 
Moreover $( \omr, W, \theta)_l$ coincides with 
$( \omr, W, \theta)_{l-1}$ on $U_{l-1}$.

By induction assume that for some $l \leq m$, we have constructed 
a real-analytic $( \omr, W, \theta)_{l-1}$ on $U_{l-1}$ which satisfies 
$(S_k)$ for all $k \leq l-1$. We now solve $S_{l}$ on some
$U_{l}$ extending $U_{l-1}$ with the initial conditions
\[ ( \omr, W, \theta)|_{t_l = 0} = ( \omr, W, \theta)_{l-1}. \]
Notice that $S_{l}$ is underdetermined,
in fact the $w^{ij}$'s with $j > l$ appear only in the last
equation. Moreover, if we ignore the last equation, the unknowns 
$w^{ij}$ and $\theta_i$ 
with $i>l$ and $j \leq l$, appear only in the second and third set of 
equations which are independent from the others.
To solve the second and third set of equations with $i > l$ we
can arbitrarily choose real analytic $w^{il}$'s extending
the initial conditions and then obtain $w^{ij}$ and $\theta_i$
by integration. Now also arbitrarily extend $w^{ij}$ with 
$i,j > l$ and set $w^{ij}$ with $j>l$ and $i \leq l$ equal
to $w^{ji}$ (this is coherent with requiring $W$ to be symmetric).
One can now see that $(S_l)$ restricted to the cases $i \leq l$
is completely determined and a unique solution exists by
Cauchy-Kowalesky on some $U_l$. 

Suppose that on $U_{l-1}$, $w^{ij} = w^{ji}$ for all 
$i$ and $j$, we now show that the same continues to be true
on $U_l$ for the solution just constructed. By construction this
is true for all $i,j$ with $i > l$ or $j>l$. So assume 
$i \leq j \leq l$.
We have
\[ \frac{\partial}{\partial t_{l}}( w^{ij}- w^{ji}) =  
   \frac{\partial w^{jl}}{\partial t_{i}}- 
             \frac{\partial w^{jl}}{\partial t_{i}} = 0, \]
where in the first equality we have used the third and fourth
equation of $(S_l)$. Now the claim follows from 
$( w^{ij}- w^{ji})|_{t_l =0} = 0$.

We now show that $( \omr, W, \theta)_{l}$ continues to satisfy
$(S_k)$ for all $k < l$ on $U_l$. First we show that the third
equation of $(S_k)$ holds. 
Assume first $j \leq k \leq i$. Then
\[ \frac{\partial}{\partial t_{l}} \left(
   \frac{\partial w^{ij}}{\partial t_{k}} - 
       \frac{\partial w^{ik}}{\partial t_{j}} \right) = 
     \frac{\partial}{\partial t_{k}} \left(
   \frac{\partial w^{ij}}{\partial t_{l}} - 
       \frac{\partial w^{il}}{\partial t_{j}} \right) = 0,
\]
where both equal signs follows from applying the third equation
of $(S_l)$. Now, assuming $j \leq i \leq k$, we have  
\begin{eqnarray*} 
 \frac{\partial}{\partial t_{l}} \left(
   \frac{\partial w^{ij}}{\partial t_{k}} - 
       \frac{\partial w^{ik}}{\partial t_{j}} \right) & = &
   \frac{\partial^2 w^{ij}}{\partial t_{k} \partial t_{l}} - 
       \frac{\partial^2 w^{kl}}{\partial t_{j}\partial t_{i}} = \\
\ & = &  \frac{\partial^2 w^{ij}}{\partial t_{k} \partial t_{l}} - 
       \frac{\partial^2 w^{lk}}{\partial t_{j}\partial t_{i}} = \\
 \ & = &   \frac{\partial}{\partial t_{k}} \left(
   \frac{\partial w^{ij}}{\partial t_{l}} - 
       \frac{\partial w^{li}}{\partial t_{j}} \right) = 0,
\end{eqnarray*}
where in the first equality we have used the fourth equation,
in the second we used symmetry, in the third we used the previous
case and in the last the third equation. 
We now show that the fourth equation of $(S_k)$ holds. Assume
$i \leq j \leq k < l$, then
\begin{eqnarray*} 
 \frac{\partial}{\partial t_{l}} \left(
   \frac{\partial w^{ij}}{\partial t_{k}} - 
       \frac{\partial w^{jk}}{\partial t_{i}} \right) & = &
   \frac{\partial^2 w^{ij}}{\partial t_{k} \partial t_{l}} - 
       \frac{\partial^2 w^{kl}}{\partial t_{i}\partial t_{j}} = \\
\ & = &  \frac{\partial^2 w^{ij}}{\partial t_{k} \partial t_{l}} - 
       \frac{\partial^2 w^{lk}}{\partial t_{i}\partial t_{j}} = \\
 \ & = &   \frac{\partial}{\partial t_{k}} \left(
   \frac{\partial w^{ij}}{\partial t_{l}} - 
       \frac{\partial w^{lj}}{\partial t_{i}} \right) = 0,
\end{eqnarray*}
where the first equality follows from the fourth equation of $(S_l)$,
the second from symmetry, the third from the third equation
of $(S_k)$ and the last from symmetry and the fourth equation
of $(S_l)$. This proves that the fourth equation of $(S_k)$
continues to hold on $U_l$. Notice also that if $i \leq k < l$,
we have
\begin{equation} \label{quatre}
  \frac{\partial w^{ik}}{\partial t_{l}} = 
     \frac{\partial w^{kl}}{\partial t_{i}} = 
       \frac{\partial w^{lk}}{\partial t_{i}} =
        \frac{\partial w^{li}}{\partial t_{k}} =
          \frac{\partial w^{il}}{\partial t_{k}}.
\end{equation}
We now prove that the second equation of $(S_k)$ continues
to hold. We have
\begin{eqnarray*}
  \frac{\partial}{\partial t_{l}} 
    \left( \left[ \frac{\partial \, \theta_i}{\partial t_{k}} 
          \right]_{(0,1)} 
                 + \overline{\partial} \, w^{ik} \right) & = &
  \left[ \frac{\partial^2 \, \theta_i}{\partial t_{k}\partial t_{l}} 
     \right]_{(0,1)} 
    + \overline{\partial} \, 
        \left(\frac{\partial w^{ik}}{\partial t_{l}} \right) = \\
\ & = &  \frac{\partial}{\partial t_{k}} 
    \left( \left[ \frac{\partial \, \theta_i}{\partial t_{l}} 
          \right]_{(0,1)} 
                 + \overline{\partial} \, w^{il} \right) = 0,
\end{eqnarray*}
where in the second equality we used the third equation of
$(S_l)$ or (\ref{quatre}) and
in the last we used the second equation of $(S_l)$. Notice
also that by symmetry and by the fact that $\theta$ has
imaginary values we have
\[ \frac{\partial \, \theta_i}{\partial t_{k}} = 
            \frac{\partial \, \theta_k}{\partial t_{i}}. \]
Using this fact it is easy to show that also the first equation
of $(S_k)$ continues to hold. We leave it to the reader. 
This concludes the proof that on some $U = U_m$ there exists a
$( \omr, W, \theta) = ( \omr, W, \theta)_m$ which satisfies
$(S_l)$ for all $l$ and such that $W$ is symmetric. Together
with (\ref{quatre}), this also proves that $( \omr, W, \theta)$
is a solution to the system in Proposition~\ref{t2:forms}
and satisfies the required initial conditions.
\end{proof}

\begin{rmk}{}
It is interesting that in the case of the circle bundle, the
Calabi-Yau $X$ constructed by this theorem is unique. This means 
that the geometry of $X$ is completely determined by the initial
data on $\pi: E \rightarrow N$. 
\end{rmk}

\section{From $T$ actions to sL submanifolds} \label{s1:exmpls}
We apply Theorems~\ref{cy:exist}
and we combine it with Theorem~\ref{s1:gold} to construct new 
families of compact sL submanifolds. All the data is assumed to
be real-analytic.

\begin{observation} {Example 1.}
Let $N$ be the complex torus $\complex{n}/ \Lambda$,
where $\Lambda$ is a $2n$-dimensional lattice generated 
by vectors $e_{1}, \ldots, e_{n},f_{1}, \ldots, f_{n}$,
where $e_{1}, \ldots, e_{n}$ is the standard basis
of $\complex{n}$ and $f_{1}, \ldots, f_{n}$ are arbitrary. 
Choose holomorphic line-bundles $L_1, \ldots, L_m$ 
with $U(1)$-connections $\theta_1, \ldots, \theta_m$. Let
$U_1, \ldots, U_m$ be the associated $U(1)$-bundles respectively.
Let $E = U_1 \oplus \ldots \oplus U_m$. Suppose
\[ \omega_{0} = \frac{i}{2} 
         \sum \, h_{jk} \, dz_{j} \wedge d\overline{z}_{k} \]
is a K\"{a}hler form on $N$ whose coefficients
$h_{jk}$ are all real valued. 
As holomorphic $n$-form we take the standard one
\[ \Omega_{0} = dz_{1} \wedge \ldots \wedge dz_{n}. \]
Define $\theta(0)=(\theta_1, \ldots, \theta_m)$ and choose a matrix 
$W_0$ satisfying the conditions of Theorem~\ref{cy:exist}. 
It then follows that we have a CY structure on $X = E \times U$ 
for some $U \subseteq \reals{m}$ with the properties of 
Theorem~\ref{cy:exist}. 
Now, if
$\Img \Lambda \subseteq \reals{n}$ is the lattice
spanned by $\Img f_{1}, \ldots, \Img f_{n}$, let 
$B = \reals{n}/ \Img \Lambda$. From our choice of
$\Omega_{0}$ and $\omega_{0}$ it follows that the
map $f: N \rightarrow B$ given by
\[ f(z_{1}, \ldots, z_{n}) = (\Img z_{1}, \ldots, \Img z_{n}) \]
is an $n$-dimensional sL fibration on $(N, \Omega(0), \omega(0))$. 
Therefore, using Theorem~\ref{s1:gold}, we can lift it to a family of 
sL $n+m$-tori of $X$.

Take for example $m=1$ and $U_1 = N \times S^1$ with $\theta_1 = i ds$, 
where $s$ is the coordinate on $S^1$. In this case, $W_0$ is a 
function uniquely determined by $\omega_0$ and $\Omega_0$ via 
(\ref{omo:Omo}). 
The data $(U_1, \Omega_0, \omega_{0}, \theta_0)$ 
can be interpreted as an $n$-parameter family of Riemannian
$n+1$-tori $T_{y}= S^{1} \times f^{-1}(y)$, where $y \in B$. 
The metric on $T_{y}$ is given by $g_{y} = 
\det( h^{jk}) \, ds^{2} + h_{jk} \, dx_{j}dx_{k}$, where the $x_{k}$'s 
are the coordinates on $f^{-1}(y)$. Theorem~\ref{cy:exist} embeds this
family isometrically into $X$.
We obtain different metrics by choosing more general $\theta_{0}$'s. 
\end{observation}

\begin{observation} {Example 2.}  The following example
is taken from Bryant~\cite{bry:realslice}.
Let $c: \complex{n} \rightarrow \complex{n}$ be the
standard conjugation. The map $c$ obviously descends
to a map of $\projc{n-1}$ (which we still call $c$).
The fixed locus of $c$ is $\projr{n-1} \subseteq \projc{n-1}$.
Let $N$ be a smooth algebraic hyper-surface of $\projc{n-1}$
with trivial canonical bundle. 
Assume $N$ to be invariant under $c$, for example
$N$ can be originating from a homogeneous polynomial of degree $n$
in $\complex{n}$ with real coefficients.
Let $Y$ be the fixed locus of $c$ restricted to $N$, i.e.
$Y = N \cap \projr{n-1}$. 
On $N$ we have a natural K\"{a}hler form $\omega_{0}$, 
namely the restriction of the Fubini-Study metric,
with respect to which $Y$ is Lagrangian. It is also 
possible to chose a nowhere vanishing holomorphic
$(n-2)$-form $\Omega_{0}$ on $N$ such that 
$c^{\ast} \Omega_{0} = \overline{\Omega}_{0}$. Given such a form
we have $\Img \Omega_{0|Y} = 0$. 
Thus $Y$ is sL with respect to the ACY structure 
$(N, \Omega_{0}, \omega_{0})$. 
We can now take $E = U_1 \oplus \ldots \oplus U_m$, where 
the $U_j$'s are holomorphic $U(1)$-bundles on $N$ with 
$U(1)$-connections $\theta_1, \ldots, \theta_m$.
We apply the corollary to construct a Calabi-Yau metric
on $X=E \times U$ for some $U \rightarrow \reals{m}$. 
Theorem~\ref{s1:gold}
then allows us to lift $Y$ to an $(n+m-2)$-dimensional sL 
submanifold.
It was shown by Bryant~\cite{bry:realslice}
that in the cases $n=3,4$ and $5$, $Y$ can also be
a real torus. 
\end{observation}

\observation{Example 3} This is a particular case of Example 1.
Let $N = \complex{\ast}$ and
$\pi: E \rightarrow N$ be a $T^2$ principal bundle
with a connection $\Theta_0$. We let
$\Omr = du$.
Given a $(2 \times 2)$ positive definite matrix of functions
$W_0$ on $N$, define 
\[
 \omega_{0} = \frac{i}{2} \det W^{-1}_{0} \, \Omega_{0} \wedge \Ombar_{0}.
\]

Then Theorem~\ref{cy:exist} gives us
a $T^2$-symmetric CY structure on 
$X = E \times (-\epsilon, \epsilon)^2$, for some $\epsilon$. 

Now consider $f: N \rightarrow \reals{}$ given 
by $u \mapsto \log |u|$.
Clearly $f$ is sL with respect to the 
ACY structure $(du, \omega_0)$. Lifting it, we obtain that
$E \times \{0 \} \times \{0 \}$ is a one parameter
family of sL manifolds inside $X$. This example therefore
could have been constructed also using Theorem~\ref{onep:sl}.
In the case $E = T^2 \times N$, we will study this 
example in more detail in the next section.  
\section{SL fibrations with $T^{2}$ symmetry} \label{sl:t2}
Let $B \subseteq \reals{3}$ 
be open and let $u=x+iy$ be the complex coordinate on $\complex{}$.
Define
\begin{equation} \label{x:defin}
   \tilde{X} = \{(\phi_1, \phi_2, u, t_1, t_2) \in 
      S^{1} \times S^{1} \times \complex{} \times \reals{} \times \reals{} \, | 
                            \,  (t_1, t_2, \Img u) \in B \} 
\end{equation} 
Let $\tilde{N} = \pi_{\complex{}}(\tilde{X})$, where $\pi_{\complex{}}$ is the
projection on the $\complex{}$ component. 
We assume that on $\tilde{X}$ we have a CY structure
such that the obvious $T^2$ action on $(\phi_{1}, \phi_{2})$ is 
structure preserving, the coordinates $t_j$ are
Hamiltonians for $\eta_j = \partial / \partial \phi_j$ and
$\tilde{N}$ is the reduced space such that $\Omr = du$. The
CY structure is then described by Proposition~\ref{t2:forms}.
We recall the main data involved. We have a $(2 \times 2)$ matrix 
$W$ and a connection $\Theta = (\theta_1, \theta_2)$, which we
write as
\[ \theta_{j} = \alpha_{j} \, du-\overline{\alpha}_{j} \, d \overline{u} 
              + i\,d \phi_{j}. \]
We then have the reduced K\"{a}hler form on $N$ given by
\[
 \omr = \frac{i}{2} \det W^{-1} \, du \wedge d \ubar.
\]
This data satisfies the equations of Proposition~\ref{t2:forms}.
The forms $\omega$ and $\Omega$ on the total space $\tilde{X}$ are
\begin{eqnarray}
 \omega & = & \pi^{\ast} \omr-i \, \theta_{j} \wedge dt_{j} =  \nonumber \\
        &  = & \frac{i}{2} \left( \det W^{-1} \, du \wedge d \ubar 
               + w_{rl}( w^{lk} dt_k - \theta_{l} ) 
                          \wedge (w^{rm} dt_m + \theta_{r}) \right), \label{t2:om} \\
 \Omega & = & - ( w^{1k} dt_k - \theta_{1}) \wedge  ( w^{2k} dt_k - \theta_{2})
                \wedge du. \label{t2:Om}
\end{eqnarray}
Assume also that $W$ and $\Theta$ are periodic
in $x = \Rl u$ of period $\kappa$. Then $(\Omega, \omega)$
defines a CY structure on $X = \tilde{X}/ \kappa \, \integ{}$,
where $\kappa \integ{}$ acts by translations on $x$. 
Let $N = \tilde{N}/ \kappa \integ{}$. Observe now that the map
\begin{eqnarray}
       f:\  X&  \longrightarrow & B \label{t2:fib} \\
   (\phi_1, \phi_2, u, t_{1}, t_{2})& \longmapsto &(t_{1}, t_{2}, \Img u) 
                    \nonumber 
\end{eqnarray}
is a sL 3-torus fibration.
We will now compute: the period matrix $P$ with respect to
$\omega$, the metric on the fibres, McLean's metric on $B$, 
the volume of the fibres and the semi-flat volumes function
$\Phi$. We refer to the Introduction for the definition of these
concepts. Denote by $F_b$ a fibre over $b \in B$.

For every $b \in B$ fix a basis $\Sigma_{1}, \Sigma_{2}, \Sigma_{3}$ 
for $H_{1}(F_b, \integ{})$
to be represented by the $1$-cycles $\{ x=\phi_2 =0 \}$,
$\{ x=\phi_1 =0 \}$ and $\{ \phi_1=\phi_2 =0 \}$ respectively.
Define
\begin{eqnarray*}
& \Xi_{k}  = (\iota_{\partial/ \partial t_{k}} \omega )_{|F_b}, \ \text{when} \ 
            k=1,2 \ \text{and} \\
& \Xi_{3}  = (\iota_{\partial/ \partial y} \omega )_{|F_b}.
\end{eqnarray*}
These are harmonic $1$-forms, by McLean's theorem.

The period matrix $P = (P_{jk})$ of the fibration with respect to 
the frame $\Xi = (\Xi_1, \Xi_2, \Xi_3)$ and basis 
$\Sigma_{1}, \Sigma_{2}, \Sigma_{3}$ 
of $H_{1}(F_b, \integ{})$ is given by
\[ P_{jk} = - \int_{\Sigma_{j}}\, \Xi_{k}.\]
We have
\begin{lem} The period matrix with respect to the frame
$\Xi = (\Xi_1, \Xi_2, \Xi_3)$ and basis
$\Sigma_{1}, \Sigma_{2}, \Sigma_{3}$ of $H_1(F_b, \integ{})$
 is 
\begin{equation} \label{t2:periods}
P = \left(  \begin{array}{ccc}
    1               & 0         & \  0 \\
       \ & \ & \ \\
    0               & 1           & 0 \\
     \ & \ & \ \\
 \int_{0}^{\kappa}2\Img \alpha_1\, dx  & \ \int_{0}^{\kappa}2\Img \alpha_2 \, dx
        & \  \int_{0}^{\kappa} \det W^{-1}\, dx 
\end{array} \right).
\end{equation}
\end{lem}
The computation is straight forward and we leave it to the reader.

Notice that $P$ depends on $b \in B$. From the period matrix $P$
we obtain the period $1$-forms on $B$:
\[ \lambda_{j} = P_{j1} \,  dt_1 +  P_{j2} \, dt_2 +  P_{j3} \, dy \] 
which uniquely determine the lattice 
$\Lambda = \text{span}_{\integ{}}( \lambda_{1}, \lambda_{3}, \lambda_{3})$
inside $T^{\ast}B$.

Standard theory of Lagrangian torus fibrations 
tells us that given a Lagrangian section $\tau: B \rightarrow X$ 
of $f$, we can naturally identify $X$ with
$T^{\ast}B / \Lambda$. The identification, in our case, goes
as follows. We know that $\eta_{j}$ is the
Hamiltonian vector field corresponding to $t_{j}$. One can
check that 
\[ \zeta = \det W \left( \delvect{x} - 2 \Img \alpha_{k} \, \eta_{k} \right) \]
is the Hamiltonian vector field corresponding to $y$. 
Let $\Phi_{\eta_{j}}^{s}$ and $\Phi_{\zeta}^{s}$
denote the flows of $\eta_{j}$ and $\zeta$ respectively.
Then, the map $\Psi: T^{\ast}B / \Lambda \rightarrow X$
given by
\[ \Psi(b, s_{1} dt_{1} + s_{2} dt_{2} + s_{3} dy) =
           \Phi_{\eta_{1}}^{s_{1}} \circ \Phi_{\eta_{2}}^{s_{2}}
                    \circ \Phi_{\zeta}^{s_{3}}( \tau(b)) \]
is well defined with respect to the quotient and provides
the identification. We have the following
\begin{thm}
With respect to the frame $(\eta_{1}, \eta_{2}, \zeta)$,
the metric of the fibre $F_{b}$ has the form
\begin{equation} \label{t2:g}
g = \left(  \begin{array}{cc}
    W               &\  0   \\
       \ & \ \\
    0               &\  \det W 
\end{array} \right).
\end{equation}
With respect to the frame $\Xi = (\Xi_1, \Xi_2, \Xi_3)$,
McLean's metric $G = (\inner{ \Xi_{j}}{\Xi_{k}}_{L^2(F_b)})$
has the form
\begin{equation} \label{t2:G}
G = \int_{0}^{\kappa} g^{-1} dx.
\end{equation}
\end{thm}
\begin{proof}{Proof.} We first compute the metric on $F_b$
with respect to the frame $(\eta_{1}, \eta_2, \partial / \partial x)$
which we denote $\tilde{g} = (\tilde{g}_{jk})$.
 Obviously when $j$ and $k$ are $1$ or $2$
we have
\[ \tilde{g}_{jk} = w_{jk}. \]
Using the fact that the forms $(w^{jk} dt_{k} - \theta_{j})$
with $j=1,2$ and $du$ are of type $(1,0)$ on $X$ and formula (\ref{t2:om})
we can compute that
\begin{eqnarray*}
\tilde{g}_{j3} = \omega(\partial / \partial x, J \eta_j)
                     =   2 \Img \alpha_k \, w_{kj},
\end{eqnarray*}
when $j=1,2$, and 
\[ \tilde{g}_{33} = 
  \left| \delvect{x} \right|^{2} = \det W^{-1} +
               4 \Img \alpha_j \Img \alpha_k \, w_{jk}. \]
One can verify that, on passing from the frame 
$(\eta_{1}, \eta_2, \partial / \partial x)$ to the frame 
$(\eta_{1}, \eta_{2}, \zeta)$, one obtains that the metric has
the form (\ref{t2:g}).

We now compute McLean's metric. We easily see that
\begin{eqnarray*}
& \Xi_{k}  = -(2 \Img \alpha_k \, dx + d\phi_k) \ \text{when} \ 
            k=1,2 \ \text{and} \\
& \Xi_{3}  = - \det W^{-1} \, dx.
\end{eqnarray*}
Notice that, as a co-frame along $F_b$, $- \Xi$ is dual
to $(\eta_{1}, \eta_{2}, \zeta)$, therefore we automatically
obtain (\ref{t2:G}).
\end{proof}

Notice that $g$ can be interpreted (via the identification $\Psi$)
as a metric on the fibres of $T^{\ast}B/ \Lambda$, the same is
true for $G^{-1}$: it may be interpreted as a flat metric on 
the fibres of $T^{\ast}B/\Lambda$. 

To compute the volume of the fibres we observe that the
volume form of $F_{b}$ is $\Omega_{|F_{b}}$, by the sL property.
Thus applying it to the frame $(\eta_{1}, \eta_2, \partial / \partial x)$
and integrating over $F_b$ we obtain
\[ \Vl(F_b) = \kappa. \]
We can now compute the semi-flat volumes function $\Phi$, mentioned in
the Introduction (formula (\ref{semf:vol})) and defined by Hitchin 
\cite{hitch:msslag}

\begin{cor}
The semi-flat volumes function $\Phi$ for $T^2$ symmetric
$3$-dimensional sL fibrations is
\begin{equation} \label{sf:vol}
 \Phi =
 \frac{ \det \left( \int_{0}^{\kappa} W^{-1} \, dx \right)}
                   { \int_{0}^{\kappa} \det W^{-1} \, dx }.
\end{equation}
\end{cor}
\begin{proof}{Proof.}
A basis of harmonic forms with respect to which the period matrix
is the identity is 
\[ \Xi^{\prime} = \Xi \cdot P^{-1}. \]
With respect to this new basis, McLean's metric is
\[ G^{\prime} = (P^{-1})^{t} G P^{-1}. \]
By definition $\Phi = \det G^{\prime}$. The reader can check
that we obtain (\ref{sf:vol}).
\end{proof}
This formula generalizes a similar one obtained in Section 5 of
\cite{diego:slonepar}.
By applying Theorem~\ref{cy:exist} we see that one can find
CY metrics on $X$ with $W$ chosen arbitrarily along the line
$\{ t_1 = t_2 = 0 \}$. Therefore one obtains examples where $\Phi$
is not constant. 

\bibliographystyle{plain}

\end{document}